\documentclass[12pt]{article}

\usepackage{graphicx}
\usepackage{amsfonts}
\usepackage{amsmath}
\usepackage{amssymb}
\usepackage{enumitem}
\usepackage{mathrsfs}

\begin{document}
\setcounter{page}{1} \setcounter{section}{0}
\newtheorem{theorem}{Theorem}[section]
\newtheorem{lemma}[theorem]{Lemma}
\newtheorem{corollary}[theorem]{Corollary}
\newtheorem{proposition}[theorem]{Proposition}
\newtheorem{observation}[theorem]{Observation}
\newtheorem{definition}[theorem]{Definition}
\newtheorem{claim}{Claim}
\newtheorem{conjecture}[theorem]{Conjecture}
\newtheorem{problem}[theorem]{Problem}

\title{Average degree in graph powers}
\author{Matt DeVos \footnote{mdevos@sfu.ca. Supported in part by an NSERC Discovery Grant and a Sloan Fellowship}\\
Jessica McDonald\footnote{jessica\textunderscore mcdonald@sfu.ca. Supported by an NSERC Postdoctoral Fellowship}
\\
Diego Scheide \footnote{dscheide@sfu.ca}\\
\medskip\\
Department of Mathematics\\
Simon Fraser University\\
Burnaby, B.C., Canada V5A 1S6}
\date{December 2, 2010}

\maketitle

\bigskip

\begin{abstract}
The $k$th power of a simple graph $G$, denoted $G^k$, is the graph with vertex set $V(G)$ where two vertices are adjacent if they are within distance $k$ in $G$. We are interested in finding lower bounds on the average degree of $G^k$.  Here we prove that if $G$ is connected with minimum degree
$d \ge 2$ and $|V(G)| \ge \frac{8}{3}d$, then $G^4$ has average degree at least $\frac{7}{3}d$.  We also prove that if $G$ is a
connected $d$-regular graph on $n$ vertices with diameter at least $3k+3$, then the average degree of $G^{3k+2}$ is at least
\[(2k+1)(d+1)  - k(k+1) (d+1)^2/n - 1.\]
Both of these results are shown to be essentially best possible; the second is best possible even when $n/d$ is arbitrarily large.
\end{abstract}

\section{Introduction}

Throughout this paper we restrict our attention to finite simple connected graphs. This allows us, in particular, to refer to the average degree $a(G)$ of a graph $G$. The $k$th power of a graph $G$, denoted $G^k$, is the graph with vertex set $V(G)$ where two vertices are adjacent if they are within distance $k$ in $G$, i.e., joined by a path of length at most $k$ in $G$. It is natural to expect that $a(G^k)$ should generally be large and for this reason we are interested in finding lower bounds on this quantity.

The maximum distance between any pair of vertices in a graph $G$ is called the diameter of $G$ and denoted $diam(G)$. If $diam(G)\leq r$, then $G^k$ is a clique for all $k\geq r$, and powers of $G$ higher than $r$ do not have any additional edges. For this reason, proving good lower bounds on $a(G^k)$ often necessitates a large diameter assumption. Indeed, the problem of counting edges in $G^k$ was first considered by Hegarty \cite{He}, who proved that for a connected $d$-regular graph $G$ with diameter at least 3,
$$a(G^3) \geq (1+c)d,$$
where $c=0.087.$ The constant $c$ was improved to $1/6$ by Pokrovskiy \cite{Po}, and then to $3/4$ by DeVos and Thomass\'{e} \cite{DT}, who also weakened the assumption to minimum degree $\delta(G)$ at least $d$. The later authors provided a family of examples proving that $3/4$ is best possible for $G^3$. In contrast to this result, when $k=2$ there is no positive constant $c$ for which $a(G^k) > (1+c)a(G)$, even in the case when $G$ is connected, regular, and has a diameter constraint (see \cite{He}).

In this paper we prove the following two new essentially best-possible lower bounds on $a(G^k)$, handling the cases  $k=4$ and $k\equiv 2$ (mod 3).

\begin{theorem}\label{thm:power4} If $G$ is a connected $n$ vertex graph with $\delta(G) \ge d$ and $n \ge \frac{8}{3}d$ for an integer $d \ge 2$, then
$$a(G^4)\ge \tfrac{7}{3}d.$$
\end{theorem}

\begin{theorem}\label{thm:2mod3}
If $G$ is a connected $d$-regular graph on $n$ vertices, $k\equiv 2$ (mod 3), and $diam(G)> k $, then
\[a(G^{k})\geq \left(\tfrac{2k-1}{3}\right)(d+1)  - \tfrac{(k-2)(k+1) (d+1)^2}{9n} - 1.\]
\end{theorem}
The proof of Theorem \ref{thm:power4} comprises Section 3 of this paper and the proof of Theorem \ref{thm:2mod3} is the subject of Section 2.

\begin{figure}[htbp]
\centering
  \includegraphics[height=1.7cm]{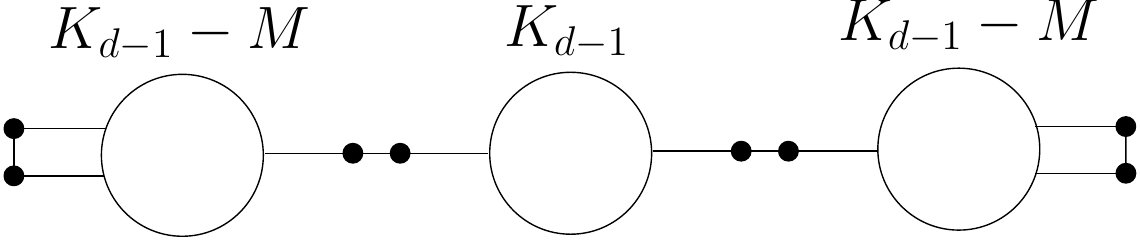}\\
   \caption{a $d$-regular graph which is extreme for Theorem \ref{thm:power4}}\label{fig:4th-pow-example}
\end{figure}

Note that the assumption $n \geq \frac{8}{3}d$ in Theorem \ref{thm:power4} could be replaced by the more restrictive $diam(G)\geq 6$. This is because a shortest path of length six contains three vertices whose neighbourhoods are completely disjoint, so $n\geq 3d>\tfrac{8}{3}t$.  In addition to proving Theorem \ref{thm:power4} in Section 3, we also present a family of examples showing the value $8/3$ cannot be further lowered. To see that the coefficient $\frac{7}{3}$ cannot be increased, consider the graph in Figure \ref{fig:4th-pow-example} (here and in later figures, each line segment represents a complete bipartite graph of the appropriate size, and each \mbox{``-M''} indicates the removal of a perfect matching). The graph is $d$-regular for
every odd $d > 1$ and a quick calculation reveals that it has $3d+4$ vertices and its fourth power has degree sum $7d^2 + 19d + 6$.

\bigskip

\begin{figure}[htbp]
\centering
  \includegraphics[height=1.7cm]{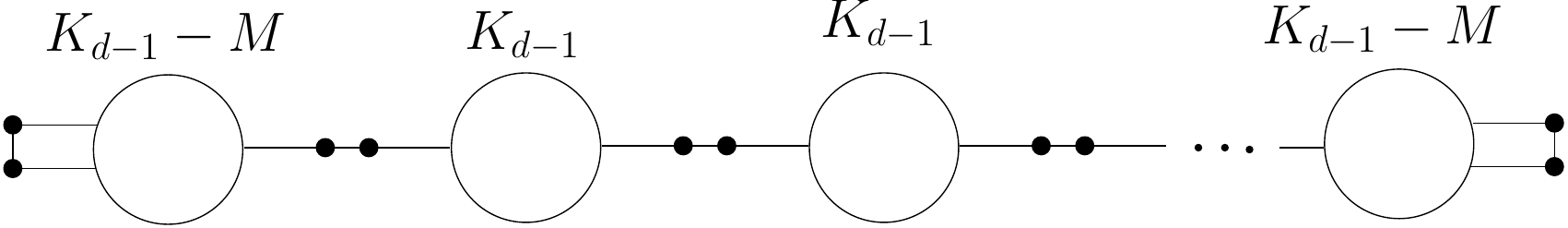}\\
   \caption{a $d$-regular graph which is extreme for Theorem \ref{thm:2mod3}}\label{fig:2mod3example}
\end{figure}

For both Theorem \ref{thm:power4} and for the $G^3$ result of DeVos and Thomass\'{e}, we know of no tight examples with arbitrarily large diameter, or with $n/d$ arbitrarily large. So it is possible that as diameter grows to infinity, better bounds could be obtained for both $G^3$ and $G^4$. The graph $G^5$, and in fact all graph powers that are 2 modulo 3, seem easier to understand. Namely, our Theorem \ref{thm:2mod3} is best possible even as $n/d$ grows to infinity. To see this, let $d >1$ be odd and consider the graph $H$ given in Figure \ref{fig:2mod3example} --- $H$ is $d$-regular graph and a generalization of Figure \ref{fig:4th-pow-example}. Similar graphs have appeared in the papers of Hegarty \cite{He} and Pokrovskiy \cite{Po}, and a straightforward calculation (which we carry out in an appendix) shows that
if $k\equiv 2$ (mod 3) and $diam(H) \geq k+1$, then
\[ a(H^{k}) \le \left(\tfrac{2k-1}{3}\right)(d+1)  - \tfrac{(k-2)(k+1) (d+1)^2}{9n} + 3. \]
This implies that Theorem \ref{thm:2mod3} cannot be improved by an additive constant greater than $4$.  On the other hand, it should be noted
that the average degree of $H^{k-1}$ is nearly that of $H^{k}$, so it seems quite possible that our theorem
could be improved by decreasing the exponent of $G$ from $k$ to $k-1$, and perhaps increasing the constant term slightly.

We can drop the assumption $diam(G)>k$ in Theorem \ref{thm:2mod3} for only a small cost. That is, when $diam(G) \leq k$ the graph $G^{k}$ is complete and we have
\begin{align*} a(G^{k}) &\ge n -1- \left( \sqrt{n} - \tfrac{1}{3\sqrt{n} }(k- \tfrac{1}{2})(d+1) \right)^2\\
	&= \left(\tfrac{2k-1}{3}\right)(d+1) - \tfrac{(k-1/2)^2(d+1)^2}{9n} - 1.
\end{align*}
Hence we get the following corollary to Theorem \ref{thm:2mod3}.

\begin{corollary}\label{cor:noDiamAsump}
If $G$ is a connected $d$-regular graph on $n$ vertices and $k\equiv 2$ (mod 3), then
\[a(G^{k})\geq \left(\tfrac{2k-1}{3}\right)(d+1) - \tfrac{(k-1/2)^2(d+1)^2}{9n} - 1.\]
\end{corollary}

We can also rewrite Corollary \ref{cor:noDiamAsump} using the parameter $diam(G)$ instead of $n$. To see this, set
$t = diam(G)$ and choose $v_0, v_1, \ldots, v_t$
to be the vertex sequence of a geodesic (shortest) path between $v_0$ and $v_t$.  Now, the neighbourhoods of
$v_0, v_3, \ldots v_{ 3\lfloor \frac{t}{3} \rfloor }$ are pairwise disjoint and it follows that $n \ge (d+1) \tfrac{diam(G)}{3}$.  Hence we get the following.

\begin{corollary}\label{cor:diamPar} If $G$ is a connected $d$-regular graph and $k\equiv 2$ (mod 3), then
\[a(G^k) \geq \left(\tfrac{2k-1}{3}\right)(d+1)\left(1-\tfrac{2k-1}{4 \, diam(G)}\right)-1.\]
\end{corollary}

For vertex transitive graphs, the bound of Theorem \ref{thm:2mod3} may be improved. Let $G$ be a finite $d$-regular vertex transitive graph and let $k < diam(G)$.  Now the degree of a vertex $x$
in $G^k$ will be $|N^1(x) \cup N^2(x) \ldots \cup N^k(x)|$ where $N^i(x) = \{ y \in V(G) : dist(x,y) = i\}$.  It is
immediate that $|N^1(x)| = d$ and since each $N^i(x)$ with $1 \le i < diam(G)$ is a vertex cut, it follows
from a theorem of Mader \cite{Mad} and Watkins \cite{Wat} (see page 40 of \cite{GR} for a proof) that $|N^i(x)| \ge \frac{2}{3}(d+1)$. This gives us the bound
\begin{equation}\label{eq:transitive}
a(G^k) \ge \tfrac{2k+1}{3}(d+1) - 1.
\end{equation}
This bound is best possible due to a graph which is constructed from a collection of $>k$ disjoint cliques of size $\frac{1}{3}(d+1)$
by placing them in a cyclic order and joining each completely to its neighbours in this ordering.  Let us note that (\ref{eq:transitive}) is still
quite close to the bound of Theorem \ref{thm:2mod3}.  Indeed, in both cases, increasing $k$ by $3$ has the effect of improving the bound
by $2(d+1)$.

This last result has consequences in additive number theory and group theory by way of Cayley graphs.  Let $\Gamma$
be a finite multiplicative group and let $A \subseteq \Gamma$ be a generating set with $1 \in A$ and with the property that
$g \in A \Rightarrow g^{-1} \in A$.  If $G$ is the Cayley graph generated by $A \setminus \{1\}$ then $G$ is a regular graph
of degree $|A| - 1$ and $G^k$ will be the Cayley graph generated by $A^k \setminus \{1\}$, so it will be regular of degree
$|A^k| - 1$.  Since Cayley graphs are vertex-transitive we may apply (\ref{eq:transitive}), which shows that whenever $A^k \neq \Gamma$
\[ |A^k| \ge \tfrac{2k+1}{3}|A| .\]
This bound is traditionally obtained in additive number theory by way of Kneser's addition theorem \cite{Kn}.

\section{Proof of Theorem \ref{thm:2mod3}}

Let $P_n$ denote a path on $n$ vertices. Our proof of Theorem \ref{thm:2mod3} relies on the number of edges $e(T^k)=|E(T^k)|$ in the $k$th power of a tree $T$.

\begin{observation}
$e(P_n^k) \ge kn - \frac{1}{2}k(k+1)$
\end{observation}

\noindent{\it Proof:}
The total degree sum in $P_n^k$ is at least $2kn - 2(1 + 2 + \cdots + k) = 2kn - k(k+1)$ so
$e(P_n^k) \ge kn - \frac{1}{2}k(k+1)$.

\begin{lemma}
\label{lem:treepower}
If $T$ is a tree on $n$ vertices, then $e(T^k) \ge kn - \frac{1}{2}k(k+1)$
\end{lemma}

\noindent{\it Proof:} Since $e(P_n^k) \ge kn - \frac{1}{2}k(k+1)$, it suffices to prove that $e(T^k) \ge e(P_n^k)$.
We prove this by induction on $\sum_{v \in V(G)} \max \{0,  \mathit{deg}(v) -2 \}$.  As a base, observe that if this sum is $0$, then
$T$ is isomorphic to $P_n$ and the result is immediate.  For the inductive step, we may then assume that
there exists a vertex of degree $\ge 3$.  Fix a root vertex $r$, and choose a vertex $v$ so that $\mathit{deg}(v) \ge 3$ and subject to this,
$v$ has maximum distance from the root. Then $T$ contains a path $P$ between two leaf vertices $u,u'$ so that $v$ is an interior vertex
of $P$, and all other interior vertices of $P$ have degree $2$ in $T$.  Let $X$ be the vertex set of this path.  Now, we modify our tree $T$
to form a new tree $U$ by deleting an edge of $P$ which is incident with $v$ and then adding the new edge $uu'$.  In the new graph $U$,
the subgraph induced by $X$ is still a path, so the number of edges in $U^k$ with both ends in $X$ is the same as that in $T^k$.  For a vertex
$w \in V(G) \setminus X$ the set of neighbours of $w$ in $V(G) \setminus X$ in the two graphs $T^k$ and $U^k$ are identical, and the number of
neighbours of $w$ in $X$ in the graph $U^k$ is at most that in $T^k$.  It follows that $e(U^k) \le e(T^k)$, and now applying the
inductive hypothesis to $U$ completes the proof.
\quad\quad$\Box$

\bigskip

We are now ready to prove Theorem \ref{thm:2mod3}, save for one very useful definition. For a graph  $G$, a vertex $v \in V(G)$, and a nonnegative integer $k$, the ball of radius $k$ around $v$ is defined to be
$B_k(v) = \{ u \in V(G) : \mathit{dist}(u,v) \le k \}$.

\begin{theorem}
If $G$ is a connected $d$-regular graph on $n$ vertices and $diam(G) > 3k+2$, then
\[a(G^{3k+2}) \ge (2k+1)(d+1) - k(k+1) (d+1)^2/n - 1.\]
\end{theorem}

\noindent{\it Proof:}
Choose a geodesic path of length $\ge 3k+3$ in $G$ and let $X_0 \subseteq V(G)=V$ consist of every third vertex of this
path.  Now, we extend $X_0$ to a set $X$ by the following procedure.  At each stage, if there exists a vertex which has distance
$3$ to $X$, then we add such a point, and otherwise we stop.  Note that $|X| \ge |X_0| \ge k+1$.  Now, construct a new graph $H$ with vertex set $X$ by the rule that $u,v \in X$ are adjacent in $H$ if they have distance $3$ in $G$.  Observe that by our construction, the graph $H$ must be connected.  We set $Z = \bigcup_{w \in X} B_1(w)$, set $Y = V \setminus Z$, and set $z = |Z|$ and $y = |Y|$ and $x = |X|$ (noting that $z = (d+1)x$).  We proceed with a sequence of claims. In what follows, $e_{3k+2}(Z,Y)$ denotes the number of edges between $Z$ and $Y$ in $G^{3k+2}$, and similarly, $e_{3k+2}(Z,Z)$ denotes the number of edges induced on $Z$ in $G^{3k+2}$.

\bigskip

\noindent{(1)} $e_{3k+2}(Z,Z) \ge (k+\frac{1}{2})(d+1)z - \frac{1}{2}z - \frac{1}{2}k(k+1)(d+1)^2$

\smallskip

First note that if $u \in X$, then $B_1(u)$ induces a clique in $G^{3k+2}$ of size $d+1$.  Next, observe that
if $u,v \in X$ are adjacent in $H^k$, then $B_1(u)$ and $B_1(v)$ will be completely joined in the graph $G^{3k+2}$, so we will have
$e_{3k+2}(B_1(u), B_1(v)) = (d+1)^2$.  Since $H$ is connected, Lemma \ref{lem:treepower} gives us
\begin{align*}
 e_{3k+2}(Z,Z)
 	&\ge \tfrac{1}{2}d(d+1)x + e(H^k)(d+1)^2		\\
	&\ge \tfrac{1}{2}d(d+1)x + (kx - \tfrac{1}{2}k(k+1))(d+1)^2 \\
	&=  (k + \tfrac{1}{2})(d+1)z - \tfrac{1}{2}z - \tfrac{1}{2}k(k+1)(d+1)^2 \\
\end{align*}
as desired.

\bigskip

\noindent{(2)} $e_{3k+2}(Z,Y) \ge k(d+1)y$

\smallskip

Let $w \in Y$ and note that by assumption, $w$ must be distance $2$ from some point $u \in X$ (were $w$ to have distance $\ge 3$ to every
point in $X$, then the set $X$ could have been augmented by adding a new point at distance 3).  Now, $|X| \ge k$ and it follows that
$\mathit{deg}_{H^{k-1}}(u) \ge k-1$.  For every point $v$ which is either equal to $u$ or a neighbour of $u$ in the graph $H^{k-1}$ we have that $w$ will be joined to $B_1(v)$ in the graph $G^{3k+2}$.  It follows from this that $|B_{3k+2}(w) \cap Z| \ge k(d+1)$ and the proof of (2)
now follows by summing this over all $w \in Y$.

\bigskip

\noindent{(3)} Every $w \in V$ satisfies $\mathit{deg}_{G^{3k+2}}(w) \ge (k+1)(d+1) - 1$

\smallskip

If $B_{3k+2}(w) = V$ then $V$ contains $\ge k+1$ disjoint balls of radius $1$ which gives the desired bound.  Otherwise, we may choose
a geodesic path of length $3k$ starting at $w$, say with vertex sequence $w = w_0, w_1, w_2, \ldots, w_{3k}$.  Now we find that
$B_{3k+2}(w)$ contains the disjoint sets $B_1(w_0)$, $B_1(w_3), \ldots, B_1(w_{3k})$ which again gives the desired bound.

\bigskip

We are now ready to complete the argument.  Below we use (1), (2), and (3) in getting to the third line.
\begin{align*}
\sum_{w \in V} \mathit{deg}_{G^{3k+2}}(w)	
	& = \sum_{w \in Z}  \mathit{deg}_{G^{3k+2}}(w) + \sum_{w \in Y}  \mathit{deg}_{G^{3k+2}}(w) \\
	& =  2 e_{3k+2}(Z,Z) + e_{3k+2}(Z,Y) +  \sum_{w \in Y}  \mathit{deg}_{G^{3k+2}}(w) \\
	&\ge (2k+1)(d+1)z - z - k(k+1)(d+1)^2  + k(d+1)y + (k+1)(d+1)y - y \\
	& =  (2k + 1)(d+1)n   - k(k+1)(d+1)^2 - n.
\end{align*}
This completes the proof.
\quad\quad$\Box$

\section{The 4th Power}

\begin{figure}[htbp]
\centering
  \includegraphics[height=2cm]{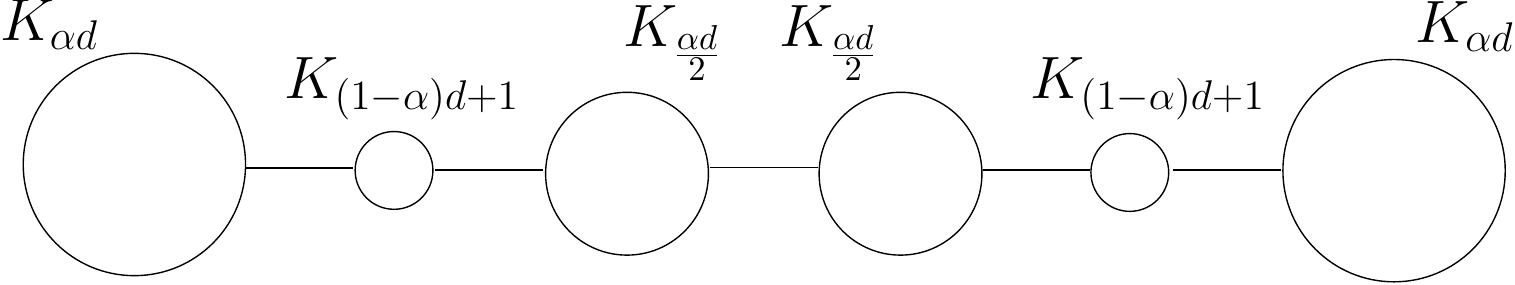}\\
   \caption{A $d$-regular graph which is extreme for Theorem \ref{thm:power4} when $\alpha <2/3$}\label{fig:83example}
\end{figure}

Before we prove Theorem \ref{thm:power4}, we give an example to show that the
value $\frac{8}{3}d$ in the theorem cannot be further lowered. To this end, consider the graph in Figure \ref{fig:83example}. This is a graph with minimum degree $d$ and $n=d(2+\alpha)+2$ vertices. We claim that if $d$ is large and $n<\frac{8}{3}d$ (and consequently $\alpha<\frac{2}{3}$), then $G^4$ has fewer than the $\frac{7}{6}nd$ edges expected by Theorem \ref{thm:power4}. Since $G^4$ is complete except for edges between vertices at distance 5, we get $e(G^4)=\frac{n(n-1)}{2}-\alpha^2d^2$.  Substituting for $n$, this gives
$$e(G^4)-\tfrac{7}{6}nd=\tfrac{1}{2}d^2\left(\alpha-\tfrac{2}{3}\right)\left(1-\alpha\right)+d\left(\tfrac{2}{3}+
\tfrac{3}{2}\alpha\right)+1.$$
This value will indeed be negative when $\alpha <2/3$, provided $d$ is chosen large enough.

\bigskip

The following lemma deals precisely with the boundary case of $\tfrac{8}{3}d\leq n\leq 3d$.

\begin{lemma}
\label{lem:smallsize}
If $G$ is a connected $n$ vertex graph with $\delta(G) \ge d$ and $\frac{8}{3}d \le n \le 3d$ for an integer $d \ge 2$, then $G^4$ has average degree $\ge \frac{7}{3}d$.
\end{lemma}

\noindent{\it Proof:} If $\mathit{diam}(G) \le 4$ then $G^4$ is complete, so it has average degree $\ge \frac{8}{3}d - 1$.  For $d \ge 3$ this is
at least $\frac{7}{3}d$, and when $d=2$ we must have $n \ge 6$ so again $G^4$ has average degree $\ge \frac{7}{3}d$.  In the remaining
case, we may choose a geodesic path with vertex sequence $v_1, v_2, \ldots, v_6$. Note that both $v_1$ and $v_6$ must have a neighbour not belonging to that path, thus we have $n\ge8$. Suppose
there exists a vertex $w$ with $\mathit{dist}(w,v_3), \mathit{dist}(w,v_4) \ge 3$. Then either
$\mathit{dist}(w,v_1) \ge 3$ or $\mathit{dist}(w,v_6) \ge 3$.  Hence  either the sets $B_1(w)$, $B_1(v_4)$, and $B_1(v_1)$ or the sets
$B_1(w)$, $B_1(v_3)$, and $B_1(v_6)$ are disjoint, which is contradictory to $n \le 3d$.  Thus, we may assume $B_2(v_3) \cup B_2(v_4) = V(G)=V$.  Now partition $V$ into the following three sets.
\begin{align*}
A &= B_2(v_3) \setminus B_2(v_4)	\\
B &= B_2(v_3) \cap B_2(v_4)	\\
C &= B_2(v_4) \setminus B_2(v_3)
\end{align*}
Set $a = |A|$ and $c = |C|$.  It is immediate from the assumption $B_2(v_3) \cup B_2(v_4) = V$ that in the graph $G^4$ the vertices in $B$ are adjacent to every other vertex and that $A$ and $C$ induce complete graphs.  Note further that $A$ and $B_1(v_3)$ and $B_1(v_6)$ are disjoint so $n - 2(d+1) \ge a$ and by a similar argument $n -2(d+1) \ge c$.  Using these observations we find
\begin{align*}
\sum_{w \in V} \mathit{deg}_{G^4}(w)
	&\ge n(n-1) - 2ac	\\
	&\ge n(n-1) - 2(n-2(d+1))^2	\\
	&= \tfrac{7}{3}nd + (n - \tfrac{8}{3}d)(3d + 6 - n) + n-8	\\
	&\ge \tfrac{7}{3}nd
\end{align*}
as desired.
\quad\quad$\Box$

\bigskip

We require one additional lemma before our proof of Theorem \ref{thm:power4}.
	
\begin{lemma}
\label{lem:bad_dist3}
Let $G$ be a connected $n$ vertex graph with $\delta(G) \ge d$ and $n> 3d$.  If there exist $u,v \in V(G)$ with $\mathit{dist}(u,v) = 3$ so
that $B_4(u) \neq V$, $B_4(v) \neq V(G)$ and $B_2(u) \cup B_2(v) = V(G)$, then $G^4$ has average degree $\ge \frac{7}{3}d$.
\end{lemma}

\noindent{\it Proof:} We define the following sets for $3 \le i \le 5$
\begin{align*}
B &= B_2(u) \cap B_2(v) \\
A_i &= \{ w \in V : \mbox{$\mathit{dist}(w,v) = i$} \} \\
C_i &= \{ w \in V : \mbox{$\mathit{dist}(w,u) = i$} \}
\end{align*}
and set $b = |B|$, $a_i = |A_i|$ and $c_i = |C_i|$.  Note that by our assumptions, these sets are
disjoint and have union equal to $V$.  Since $B_2(u) \cup B_2(v) = V(G)=V$, every point in $B$ is adjacent to every other vertex in $G^4$.
For a vertex $w \in A_3$ we have that $B_4(w)$ contains the disjoint sets $B_1(u)$, $B_1(v)$,
and $A_5$, so it will have degree $\ge 2 d + a_5$ in $G^4$.  These two observations plus a similar one for $C_3$ give us
\[ \sum_{w \in A_3 \cup B \cup C_3} \mathit{deg}_{G^4}(w) > 3 d b + 2 d(a_3 + c_3) + a_3 a_5 + c_3 c_5 \]
If $w \in A_4$, we may choose $w' \in V$ so that $\mathit{dist}(w,w') = 3$ and $\mathit{dist}(w',v) = 1$.  Now $B_4(w)$ contains the disjoint sets $A_5 \cup A_4 \cup A_3$ and $B_1(w')$ so $\mathit{deg}_{G^4}(w) \geq d + a_3 + a_4 + a_5$.  If $w \in A_5$ we may choose $w' \in V$ so that $\mathit{dist}(w,w') = 3$ and $\mathit{dist}(w',v) = 2$.  Now $B_4(w)$ contains the disjoint sets $A_5 \cup A_4$ and $B_1(w')$, so
$\mathit{deg}_{G^4}(w) \ge a_5 + a_4 + d$.  This gives us
\[ \sum_{w \in A_4 \cup A_5} \mathit{deg}_{G^4}(w) \ge (d + a_4 + a_5)(a_4 + a_5) + a_3 a_4. \]
By a similar argument we get
\[ \sum_{w \in C_4 \cup C_5} \mathit{deg}_{G^4}(w) \ge (d + c_4 + c_5)(c_4 + c_5) + c_3 c_4. \]
Now, $B_4(u) \neq V$ so there exists a point $w \in C_5$ and $B_1(w) \subseteq C_4 \cup C_5$.  Thus $c_4 + c_5 \ge d$ and $c_3 c_4 + c_3 c_5 \ge c_3 d$.  Similarly $a_3 a_4 + a_3 a_5 \ge a_3 d$.
Setting $\bar{a} = a_4 + a_5$ and $\bar{c} = c_4 + c_5$ and combining our above inequalities with these observations yields
\begin{align*}
\sum_{w \in V} \mathit{deg}_{G^4}(w) - \tfrac{7}{3} d n		
	&\ge 3 d( n - \bar{a} - \bar{c} ) + (d + \bar{a}) \bar{a} + (d + \bar{c} ) \bar{c} - \tfrac{7}{3} d n \\
	&= \tfrac{2}{3} d n  - 2 d (\bar{a} + \bar{c}) + \bar{a}^2 + \bar{c}^2	\\
	&\ge 2 d^2 - 2 d( \bar{a} + \bar{c} ) + \frac{ (\bar{a} + \bar{c} )^2 }{2}	\\
	& =  \left( \sqrt{2} d - \frac{ \bar{a} + \bar{c} }{ \sqrt{2} } \right)^2 \\
	& \ge 0
\end{align*}
which completes the proof.
\quad\quad$\Box$

\begin{theorem}
If $G=(V,E)$ is a connected $n$ vertex graph with $\delta(G) \ge d$ and $n \ge \frac{8}{3}d$ for an integer $d \ge 2$, then $G^4$ has average degree $\ge \frac{7}{3}d$.
\end{theorem}

\noindent{\it Proof:} Let $G$ be a counterexample with $n$ minimum.  Define a vertex $v$ to be \emph{good} if $\mathit{deg}_{G^4}(v) \ge 3 d$ and \emph{bad} otherwise.  Let $Z \subseteq V$ be the set of good vertices, and set $\gamma = \frac{|Z|}{n}$.  We prove the result with a sequence of claims.

\bigskip

\noindent{(1)} $n > 3d$

\smallskip

This follows from Lemma \ref{lem:smallsize}.

\bigskip

\noindent{(2)} Every $v \in V$ satisfies ${\mathit deg}_{G^4}(v) \ge 2 d$.

\smallskip

If $B_2(v) = V$ then the above inequality follows from (1).  Otherwise there exists $u \in V$ with ${\mathit dist}(u,v) = 3$ and now
$B_4(v)$ contains the disjoint sets $B_1(u)$ and $B_1(v)$ so $\mathit{deg}_{G^4}(v) \ge 2 d$.

\bigskip

\noindent{(3)} $G^4$ has average degree $\ge (2 + \gamma)  d$.

\smallskip

This is a consequence of the following calculation.
\begin{align*}
\sum_{v \in V} \mathit{deg}_{G^4}(v)
	&= \sum_{v \in Z} \mathit{deg}_{G^4}(v) + \sum_{v \in V \setminus Z} \mathit{deg}_{G^4}(v)  \\
	&\ge (\gamma n) (3 d) + (1-\gamma)n(2 d) \\
	&= (2 + \gamma) d n
\end{align*}

\bigskip

\noindent{(4)} If $u,v,v' \in V$ satisfy $\mathit{dist}(v,v') \ge 3 = \mathit{dist}(u,v) = \mathit{dist}(u,v')$ then
$u$ is good.

\smallskip

The sets $B_1(u)$, $B_1(v)$, and $B_1(v')$ are pairwise disjoint and are all contained in $B_4(u)$ so
$\mathit{deg}_{G^4}(u) \ge 3 d$ and $u$ is good.

\bigskip

\noindent{(5)} There do not exist bad vertices $u_1,u_2$ with $\mathit{dist}(u_1,u_2) = 3$.

\smallskip

If $u_1, u_2$ are bad, then since $n\geq 3d$, it follows from Lemma \ref{lem:bad_dist3} that we may assume that $V \setminus (B_2(u_1) \cup B_2(u_2)) \neq \emptyset$ and it
follows that there exists a vertex $w$ so that $\min \{ \mathit{dist}(u_1,w), \mathit{dist}(u_2,w) \} = 3$.  But then this a contradiction as (4) implies that one of
$u_1,u_2$ is good.

\bigskip

Let $X_1,\ldots,X_k$ be the vertex sets of the components of $G - Z$.

\bigskip

\noindent{(6)} Every $X_i$ induces a clique in $G^2$.

\smallskip

Let $v \in X_i$ and suppose (for a contradiction) that $X_i \not\subseteq B_2(v)$.  In this case, we may choose a vertex
$u \in X_i \setminus B_2(v)$ which is adjacent to a point in $B_2(v)$.  This gives us $\mathit{dist}(u,v) = 3$ contradicting (5).

\bigskip

We now define a relation on $\{X_1,\ldots,X_k\}$ by the rule that $X_i \sim X_j$ if $N(X_i) \cap N(X_j) \neq \emptyset$.

\bigskip

\noindent{(7)} If $X_i \sim X_j$ then $X_i \cup X_j$ is a clique in $G^2$.

\smallskip

Let $v \in X_i$ satisfy $N(v) \cap N(X_j) \neq \emptyset$ and suppose (for a contradiction) that $X_j \not\subseteq B_2(v)$.  Then
we may choose a vertex $u \in X_j \setminus B_2(v)$ which is adjacent to a point in $B_2(v)$.  But then $u$ and $v$ have
distance $3$ contradicting (5).  It follows that every point in $X_j$ is distance $2$ from $v$ and then by a similar argument
has distance $2$ to any point in $X_i$.

\bigskip

\noindent{(8)} $\sim$ is an equivalence relation.

\smallskip

It is immediate from the definitions that $\sim$ is both reflexive and symmetric.  To see that it is transitive, we suppose that
$X_i \sim X_j \sim X_k$.  If every point in $G$ is distance $\le 2$ to $X_i \cup X_j \cup X_k$ then it follows from (6) and (7) that
in the graph $G^4$ every point in $X_j$ is adjacent to every other vertex, but this contradicts (1) and the assumption that these
vertices are bad.  It follows that we may choose a vertex $w \in V$ so that $\mathit{dist}(w,X_i \cup X_j \cup X_k) = 3$.  First suppose that there exists $v \in X_j$ so that $\mathit{dist}(w,v) = 3$.  Now, choose $u \in X_i$ and $u' \in X_k$.  It follows from (7) that $B_4(v)$ contains $B_1(u) \cup B_1(u') \cup B_1(w)$ and since $v$ is bad this implies that $B_1(u) \cap B_1(u') \neq \emptyset$, so $X_i \sim X_k$.  Thus, we may
assume without loss that $\mathit{dist}(w,v) = 3$ for some $v \in X_i$ and that $\mathit{dist}(w,X_j) \ge 4$.  Now choose a vertex $u \in N(X_j) \cap N(X_k)$.  We must have $\mathit{dist}(u,w) \ge 3$ (otherwise
$\mathit{dist}(w,X_j) \le 3$) and $\mathit{dist}(v,u) \le 3$ by (7).  It follows that $B_4(v)$ contains $B_1(v) \cup B_1(w) \cup B_1(u)$.  Since
$v$ is bad it must be that $\mathit{dist}(u,v) \le 2$.  But then we have that $\mathit{dist}(v,X_k) \le 3$.  If there is a point in $X_k$ which is distance 3 from $v$ we get a contradiction to (5), so we must have $\mathit{dist}(v,X_k) = 2$ which implies $X_i \sim X_k$ as desired.

\bigskip

We now define $\{ Y_1, \ldots, Y_{\ell} \}$ to be the unions of the equivalence classes of $\sim$.  Note that by (6) and (7) every $Y_i$ induces a clique in $G^2$.

\bigskip

\noindent{(9)} If $1 \le i < j \le \ell$ then $N(Y_i)$ and $N(Y_j)$ are disjoint, and there are no edges between them.

\smallskip

It is immediate that $N(Y_i)$ and $N(Y_j)$ are disjoint.  Were there to be an edge between $u \in N(Y_i)$ and $v \in N(Y_j)$
then $u,v \in Z$ and we may choose $u' \in Y_i$ and $v' \in Y_j$ so that $uu', vv' \in E$.  But then we have that $u'$ and $v'$ have
distance $3$ which contradicts (5).

\bigskip

Let $Z^* = \{ u \in Z : N(u) \subseteq Z \}$.

\bigskip

\noindent{(10)} Every $v \in Z$ satisfies $|B_4(v) \cap Z| \ge d + 1$.

\smallskip

We claim that $B_3(v) \cap Z^* \neq \emptyset$ which immediately yields (10).  To show this claim, let us suppose
(for a contradiction) that it is false and choose $1 \le i \le \ell$ so that $v$ has a neighbour, say $u$, in $Y_i$.
It now follows from (9) that $B_3(v) \subseteq Y_i \cup N(Y_i)$.  However, $Y_i \subseteq B_2(u)$ so
$Y_i \cup N(Y_i) \subseteq B_3(u)$ giving us $B_3(v) \subseteq B_3(u)$.  But this contradicts the assumptions that $u$ is
bad but $v$ is good.

\bigskip

For every $1 \le i \le \ell$ and positive integer $t$ let $Y_i^t = \{ v \in Y_i : \mathit{dist}(v, Z^*) = t \}$.

\bigskip

\noindent{(11)} For $1 \le i \le \ell$ we have $Y_i^t = \emptyset$ whenever $t=1$ or $t >4$.

\smallskip

It is immediate from the definitions that $Y_i^1 = \emptyset$.  It follows from (9) that there exists a vertex $u \in N(Y_i)$ with
a neighbour in $Z^*$.  Choose $v \in Y_i$ adjacent to $u$.  Since $Y_i$ induces a clique in $G^2$ every point in $Y_i$
must have distance at most four to $Z^*$ as desired.

\bigskip

\noindent{(12)} Every $v \in Y_i^2$ satisfies ${\mathit deg}_{G^4}(v) \ge 2d + |Y_i^4|$.

\smallskip

Choose a vertex $u \in Z^*$ with $\mathit{dist}(u,v) = 2$. There must be a vertex $u'\in B_2(u)$ with $\mathit{dist}(u',v) = 3$ (otherwise $B_2(u)\subseteq B_2(v)$, contradicting that $u$ is good and $v$ is bad). Now $B_4(v)$ contains $B_1(u') \cup B_1(v) \cup Y_i^4$, and we claim that these sets are disjoint (which obviously yields (12)).
It is immediate that $B_1(u') \cap B_1(v) = \emptyset$. No point in $Y_i^4$ could be adjacent to $v$ or $u'$ since $v$ and $u'$ are distance $\le 2$ from $Z^*$, so $B_1(v) \cap Y_i^4 = B_1(u') \cap Y_i^4 = \emptyset$.

\bigskip

\noindent{(13)} For every $1 \le i \le \ell$ and $v \in Y_i \setminus Y_i^4$ we have $|B_4(v) \cap Z| \ge d$.

\smallskip

By our definitions $v$ must have distance $\le 3$ to some vertex $u \in Z^*$ but then $B_1(u) \subseteq B_4(v)$
and $B_1(u)$ is a subset of $Z$ with size $\ge d$.

\bigskip

\noindent{(14)} For every $1 \le i \le \ell$ and $v \in Y_i^4$ we have $|B_4(v) \cap Z| \ge d - |Y_i^2|$.

\smallskip

By our definitions $v$ is distance 3 to a point $u\in Z$ which has a neighbour in $Z^*$. For $j\ne i$ we have $B_1(u)\cap Y_j=\emptyset$, otherwise (9) would imply that a vertex $u'$ satisfying $\mathit{dist}(u,u') = 1$ and $\mathit{dist}(u',v) = 2$ belongs to $Z^*$, contradicting the fact that $v$ is distance 4 from $Z^*$. It follows that $B_1(u) \subseteq Z \cup Y_i^2$ and thus $B_1(u) \cap Z$ is a set of size $\ge d - |Y_i^2|$ which is contained in $B_4(v) \cap Z$.

\bigskip

Let $Y = V \setminus Z$, and for any pair of disjoint sets $S,T \subseteq V$ and positive integer $k$ we let $e_k(S,T)$ denote the number of edges between the sets $S$ and $T$ in the graph $G^k$.

\bigskip

\noindent{(15)} The average degree of $G^4$ is at least $(3 - 2\gamma) d$.

\smallskip

This is a consequence of the following equation (here we use (2), (10) and (12) in getting to the third line and (13) and (14) in
getting to the fifth line).
\begin{align*}
\sum_{v \in V} \mathit{deg}_{G^4}(v)
	&=		\sum_{v \in Z} \mathit{deg}_{G^4}(v) + \sum_{u \in Y} \mathit{deg}_{G^4}(u)	\\
	&= 		\sum_{v \in Z} e_4(v,Z \setminus \{v\}) + e_4(Z,Y) +  \sum_{u \in Y} \mathit{deg}_{G^4}(u)	\\
	&\ge		d |Z| + e_4(Z,Y)  + \sum_{i=1}^{\ell} |Y_i^2| |Y_i^4| + 2 d |Y|	\\
	&= 		d \gamma n + \sum_{v \in Y} |B_4(v) \cap Z| + \sum_{i=1}^{\ell} |Y_i^2| |Y_i^4| + 2 (1-\gamma) d n \\
	&\ge		d \gamma n + \sum_{i=1}^{\ell} |Y_i \setminus Y_i^4| d + \sum_{i=1}^{\ell} |Y_i^4| ( d - |Y_i^2| ) + \sum_{i=1}^{\ell} |Y_i^2| |Y_i^4| + 2 (1-\gamma) d n \\
	&=		d \gamma n + 3 (1 - \gamma) d n	\\
	&=		(3 - 2 \gamma) d n
\end{align*}

We can now complete the proof.  By taking a convex combination of the bounds in (3) and (15) we have that the average degree of
$G^4$ must be at least $\frac{2}{3} (2 + \gamma) d + \frac{1}{3} (3 - 2\gamma) d = \frac{7}{3} d$ thus giving us a final contradiction.
\quad\quad$\Box$

\section*{Appendix}

Here we carry out the calculation claimed in the introduction giving an upper bound on the average degree of the $k$th power (when $k\equiv 2$ (mod 3)) of the graphs appearing in Figure \ref{fig:2mod3example}.  For every odd integer $d > 1$ and every positive integer $t$ we shall define such a graph $H_t$ (where  $t+1$ is the number of large circles in the picture).  The vertex set has a partition as $\{X_{-1}, X_0, X_1, \ldots, X_{3t+1} \}$.  The edges are defined
as follows.  For $-1 \le i \le 3t$ there
is a complete bipartite graph between $X_i$ and $X_{i+1}$.  The sets $X_{-1}$ and $X_{3t+1}$ induce $K_2$, the sets
$X_0$ and $X_{3t}$ induce $K_{d-1}$ minus a perfect matching, the set $X_{3i}$ induces $K_{d-1}$ for $1 \le i \le t-1$, and $X_i$ is a single point when $3$ does not divide $i$ and $0 \le i \le 3t$.

For $0 \le i \le t$ let $Y_i$ be the union of $X_{3i}$ together with one vertex from $X_{3i-1}$ and one vertex from $X_{3i+1}$, set
$Y = \bigcup_{i=0}^t Y_i$ and set $\{u,u'\} = V(H_t) \setminus Y$.  Now, let $k$ be an integer with $k\equiv 2$ (mod 3) and $k < 3t$.  To assist in counting the edges in $H_t^{k}$ we construct the alternate graph $H_t'$ with vertex partition $\{Y_0,Y_1,\ldots,Y_t\}$ and with edges given by the rule that each $Y_i$ induces a clique, and
$Y_i$ and $Y_j$ are completely joined for $0 \le i < j \le t$ if $j-i \le \frac{k-2}{3}$.  We then find
\[ e(H_t') = (t+1)\tfrac{d(d+1)}{2} + \left( \frac{k-2}{3}(t+1) - \frac{1}{2}\frac{k-2}{3}\frac{k+1}{3} \right)(d+1)^2. \]
A vertex in $X_{3i}$ will have degree at most 2 larger in $H_t^{k}$ than in $H_t'$ while a vertex in $X_i \setminus \{u,u'\}$ with
$i$ not a multiple of $3$ will have degree at most $d+1$ larger in $H_t^{k}$ than in $H_t'$.
Using the fact that $u$ and $u'$ are vertices of minimum degree in $H_t^{k}$ this gives us
\begin{align*}
a(H_t^{k})
	&\le	\frac{1}{(t+1)(d+1)} \sum_{v \in Y} \mathit{deg}_{H_t^{k}}(v)	\\
	&\le	\frac{1}{(t+1)(d+1)} \left( \sum_{v \in Y} \mathit{deg}_{H_t'}(v) + (t+1)(d-1)2 + 2(t+1)(d+1)  \right) \\
	&\le  \tfrac{2}{(t+1)(d+1)} e(H_t') + 4	\\
	&=  (\tfrac{2k-1}{3})(d+1) - \tfrac{(k-2)(k+1)(d+1)}{9(t+1)} + 3\\
	&\le  (\tfrac{2k-1}{3})(d+1) - \tfrac{(k-2)(k+1)(d+1)^2}{9n} + 3
\end{align*}
as claimed in our introduction.


\begin{thebibliography}{00}

\bibitem{Ca} A. L. Cauchy. Recherches sur les nombres, J. \'{E}cole polytech. 9 (1813) 99-116.

\bibitem{Da} H. Davenport. On the addition of residue classes, J. London Math. Soc. 10 (1935) 30-32.

\bibitem{DT} M. DeVos and S. Thomass\'{e}. Edge Growth in Graph Cubes, preprint, 2010.

\bibitem{GR} G. Godsil and G. Royle. Algebraic graph theory, Springer, New York, 2001.

\bibitem{He} P. Hegarty. A Cauchy-Davenport type result for arbitrary regular graphs, preprint, 2009.

\bibitem{Kn} M. Kneser.  Abschätzung der asymptotischen Dichte von Summenmengen, Math. Z. 58, (1953). 459–484.

\bibitem{Mad} W. Mader. \"{U}ber den Zusammenhang symmetrischer Graphen, Arch. Math 22 (1971) 333-336.

\bibitem{Po} A. Pokrovskiy. Growth of graph powers, preprint, 2010.

\bibitem{Su} B. D. Sullivan. A summary of results and problems related to the Caccetta-H\"{a}ggkvist Conjecture, AIM Preprint 2006-13, 2006.

\bibitem{Wat} M.E. Watkins. Connectivity of transitive graphs, J. Comb. Theory 8 (1970) 23-29.


\end{thebibliography}
\end{document}